\newcommand{\rem}[1]{}
\documentclass{amsart}
\usepackage{amsfonts,amssymb,amsmath,amsthm,mathrsfs}
\usepackage{url}
\usepackage[dvips]{epsfig}
\newtheorem{thrm}{Theorem}[section]

\newtheorem{prop}[thrm]{Proposition}
\newtheorem{cor}[thrm]{Corollary}
\newtheorem{remark}[thrm]{Remark}
\theoremstyle{definition}
\newtheorem{definition}[thrm]{Definition}

\begin{document}

\author[C.~A.~Mantica and L.~G.~Molinari]
{Carlo~Alberto~Mantica and Luca~Guido~Molinari}
\address{C.~A.~Mantica: I.I.S. Lagrange, Via L. Modignani 65, 
20161, Milano, Italy -- L.~Molinari (corresponding author): Physics Department,
Universit\'a degli Studi di Milano and I.N.F.N. sez. Milano,
Via Celoria 16, 20133 Milano, Italy.}
\email{carloalberto.mantica@libero.it, luca.molinari@unimi.it}
\title[Riemann compatible tensors]
{Riemann compatible tensors}
\begin{abstract} Derdzinski and Shen's theorem on the restrictions 
posed by a Codazzi tensor on the Riemann tensor holds more generally when 
a Riemann-compatible tensor exists. Several properties are shown to 
remain valid in this broader setting.
Riemann compatibility is equivalent to the Bianchi identity 
of the new ``Codazzi deviation tensor'', with a geometric significance.  
The general properties are studied, with their implications on Pontryagin 
forms. Examples are given of manifolds with Riemann-compatible tensors,
in particular those generated by geodesic mapping. 
Compatibility is extended to
generalized curvature tensors, with an application to Weyl's tensor
and general relativity. 
\end{abstract}
\date{23 march 2012}
\maketitle
\section{Introduction}
The Riemann tensor $R_{ijk}{}^m$ and its contractions, %the Ricci tensor 
$R_{kl}=R_{kml}{}^m$ and %the scalar curvature 
$R=g^{kl}R_{kl}$,
are the fundamental tensors to describe the local structure of
a Riemannian manifold $(\mathscr M_n,g)$ of dimension $n$.
In a remarkable theorem \cite{[7],Besse} Derdzinski and Shen showed that
the existence of a non trivial Codazzi tensor poses strong constraints
on the structure of the Riemann tensor. Because of their geometric relevance, 
Codazzi tensors have been studied by several authors, as Berger and Ebin 
\cite{Berger}, Bourguignon \cite{Bou}, Derdzinski \cite{[5],[6]}, Derdzinski 
and Shen \cite{[7]}, Ferus \cite{Ferus}, Simon \cite{Simon}; a compendium of 
results is found in Besse's book \cite{Besse}.
Recently, we showed \cite{manticaDS} that the Codazzi differential condition 
\begin{equation}\label{Codazzi}
 \nabla_i b_{jk}-\nabla_j b_{ik}=0
\end{equation} 
is sufficient for the theorem to hold, and can be replaced by the 
more general notion of {\em Riemann-compatibility}, which is instead 
algebraic:
\begin{definition}
A symmetric tensor $b_{ij}$ is Riemann compatible ($R$-compatible) if:
\begin{equation}
b_{im} R_{jkl}{}^m+b_{jm} R_{kil}{}^m+b_{km} R_{ijl}{}^m=0.\label{RC}
\end{equation}
\end{definition}
With this requirement, we proved the following extension of 
Derdzinski-Shen's theorem: 
\begin{thrm}\cite{manticaDS}\label{extDS}
Suppose that a symmetric $R$-compatible tensor $b_{ij}$ exists.
Then, if $X$, $Y$ and $Z$ are three eigenvectors of the matrix $b_r{}^s$
at a point of the manifold, with eigenvalues $\lambda $, $\mu $ and $\nu $,
it is $R_{ijkl}\,X^i \, Y^j\, Z^k \, =0 $
%\begin{equation}
% R_{ijkl}\,X^i \, Y^j\, Z^k \, =0 \label{eq2.10}
%\end{equation}
provided that both $\lambda $ and $\mu $ are different from $\nu$. 
\end{thrm}

The concept of compatibility allows for a further extension of the theorem, 
where the Riemann tensor $R$ is replaced by a generalized 
curvature tensor $K$, and $b$ is required to be $K$-compatible
\cite{manticaDS}.\\

This paper studies the properties of Riemann compatibility, and its 
implications on the geometry of the manifold. In section 2 
$R$-compatibility is shown to be equivalent to the Bianchi identity of 
a new tensor, the {\em Codazzi deviation}. In section 3 the irreducible 
components of the covariant derivative of a symmetric tensor are classified in
a simple manner, based on the decomposition into traceless terms. This
is of guidance in the study of different structures suited for 
$R$-compatibility. The general properties of Riemann compatibility are 
presented in section 4. 
In section 5 several properties of manifolds in presence of 
a Riemann compatible tensor that were obtained by Derdzinsky-Shen and 
Bourguignon for manifolds with a Codazzi tensor, are recovered. 
In particular, it is shown that $R$-compatibility implies pureness, 
a property of the Riemann tensor introduced by Maillot that implies 
the vanishing of Pontryagin forms. 
Manifolds that display $R$-compatible tensors are presented in section 6; 
interesting examples are generated by geodesic mappings, that induce 
metric tensors that are $R$-compatible. Finally, in section 7, 
$K$-tensors and $K$-compatibility are presented, with applications 
to the standard curvature tensors. In the end, an application to general 
relativity is mentioned, that will be discussed fully elsewhere.

\section{The Codazzi deviation tensor and R-compatibility}
Since Codazzi tensors are Riemann compatible, for a non Codazzi
differentiable simmetric tensor field $b$ it is useful to define its 
deviation from the Codazzi condition. This tensor solves an unexpected relation 
that generalizes Lovelock's identity for the Riemann tensor, and
shows that Riemann compatibility is a condition for closedness of certain 
2-forms. 
\begin{definition} 
The {\em Codazzi deviation} of a symmetric tensor $b_{kl}$ is 
\begin{equation}
\mathscr C_{jkl}=:\nabla_j b_{kl}-\nabla_k b_{jl}
\end{equation}
\end{definition}
\noindent
Simple properties are: $\mathscr C_{jkl}=-\mathscr C_{kjl}$ and 
$\mathscr C_{jkl}+\mathscr C_{klj}+\mathscr C_{ljk} =0$.\\

The following identity holds in general, and relates the Bianchi differential
combination for $\mathscr C$ to the Riemann compatibility of $b$: 
\begin{prop}
\begin{equation}
\nabla_i\mathscr C_{jkl}+\nabla_j\mathscr C_{kil}+\nabla_k\mathscr C_{ijl} =
b_{im} R_{jkl}{}^m+b_{jm} R_{kil}{}^m+b_{km} R_{ijl}{}^m\label{CCCcomp}
\end{equation}
\begin{proof}
\begin{align*}
&\nabla_i\mathscr C_{jkl}+\nabla_j\mathscr C_{kil}+\nabla_k\mathscr C_{ijl} 
=[\nabla_i,\nabla_j]b_{kl}+[\nabla_k,\nabla_i]b_{jl}+[\nabla_j,\nabla_k]b_{il}\\
&= b_{ml}(R_{ijk}{}^m+R_{kij}{}^m+R_{jki}{}^m)+
b_{im} R_{jkl}{}^m+b_{jm} R_{kil}{}^m+b_{km} R_{ijl}{}^m
\end{align*}
the first term vanishes by the first Bianchi identity.
\end{proof}
\end{prop}

\begin{remark}\label{rem111}
The identity holds true if $b_{ij}$ is replaced by 
$b'_{ij}=b_{ij}+\chi a_{ij}$, where $a_{ij}$ is a Codazzi tensor and $\chi $ 
a scalar field. Then: $\mathscr C'_{jkl}=\mathscr C_{jkl} -
(a_{kl}\nabla_j -a_{jl} \nabla_k)\chi $.
\end{remark}

The deviation tensor is associated to the 2-form $\mathscr C_l
=\frac{1}{2}\mathscr C_{jkl}dx^j\wedge dx^k$. The closedness condition 
$0=D\mathscr C_l =\frac{1}{2}\nabla_i\mathscr C_{jkl}dx^i\wedge dx^j\wedge dx^k$
($D$ is the exterior covariant derivative)
is the second Bianchi identity for the Codazzi deviation: 
$\nabla_i\mathscr C_{jkl}+\nabla_j\mathscr C_{kil}+
\nabla_k\mathscr C_{ijl} =0 $. This gives a 
geometric picture of Riemann compatibility:
\begin{thrm}\label{pi}
$b_{ij}$ is Riemann compatible if and only if $D\mathscr C_l =0 $. 
\end{thrm}
\begin{remark} The Codazzi deviation of the Ricci tensor is, by the 
contracted second Bianchi identity:
$ \mathscr C_{jkl} =: \nabla_jR_{kl}-\nabla_k R_{jl} = -\nabla_m R_{jkl}{}^m $.
For the Ricci tensor the identity \eqref{CCCcomp} identifies with
Lovelock's identity \cite{Lovelock} for the Riemann tensor:
\begin{eqnarray}
&&\nabla_i\nabla_m R_{jkl}{}^m +\nabla_j\nabla_m R_{kil}{}^m +\nabla_k\nabla_m 
R_{ijl}{}^m \label{lovelock}\\
&&\qquad = - R_{im} R_{jkl}{}^m - R_{jm} R_{kil}{}^m - R_{km} R_{ijl}{}^m .
\nonumber
\end{eqnarray}
\end{remark}
%il lato sinistro è la derita covariante esterna della forma 
%appropriata e che questa formula e le sue varie generalizzazioni sono
%riferiti come formula di Weitenzbok per tensori curvature like 
%( come ci ha detto di mettereDerdzinski nell Coll math)

A Veblen-like identity holds, that corresponds to \eqref{veblenb}
(For $b_{ij}=R_{ij}$ %eq.\eqref{VeblenC} 
it specializes to Veblen's identity for the 
divergence of the Riemann tensor \cite{[12]}):
\begin{prop}
\begin{align}\label{VeblenC}
\nabla_i \mathscr C_{jlk} + \nabla_j \mathscr C_{kil} + 
\nabla_k \mathscr C_{lji} + \nabla_l \mathscr C_{ikj}\\ 
= b_{im} R_{jlk}{}^ m 
+ b_{jm} R_{kil}{}^m + b_{km} R_{lji}{}^m + b_{lm} R_{ikj}{}^m\nonumber
\end{align}
\begin{proof}
Write four equations \eqref{CCCcomp} with cycled indices $i,j,k,l$ and sum up. 
Then simplify by means of the first Bianchi identity for the Riemann tensor 
and the cyclic 
identity $\mathscr C_{jkl} + \mathscr C_{klj} + \mathscr C_{ljk} = 0$.
\end{proof}
\end{prop}
%%%%%%%%%%%%%%%%%%%%%%%%%%%%%%%%%%%%%%%%%%%%%%%%%%%%%%%%%%%%%%%%
%%%%%%%%%%%%%%%%%%%%%%%%%%%%%%%%%%%%%%%%%%%%%%%%%%%%%%%%%%%%%%%%

\section{Irreducible components for $\nabla_j b_{kl}$ and $R$-compatibility}
We begin with a simple procedure to classify the $O(n)$ invariant
components of the tensor $\nabla_j b_{kl}$. They will guide us in the 
study of $R$-compatibility.\\
If $b$ is the Ricci tensor, this simple construction reproduces
the seven equations linear in $\nabla_iR_{jk}$, invariant 
for the $O(n)$ group, that are enumerated and discussed in Besse's treatise 
``Einstein Manifolds'' \cite{Besse}.\\

%Cases 1) and 6) are $R$-compatible and were discussed. Case 2) (generalized
%Sinyukov), \eqref{nablabg}, is $R$-compatible if $A_j-B_j$ is closed.
%Case 8) is $R$-compatible if $\nabla_m b^m{}_j$ is closed.

%3) is equivalent to $\nabla_m R_{ijk}{}^m=0$; 4) and 7) imply that 
%$\nabla_m C_{ijk}{}^m=0$ ($C_{ijk}{}^m$ is Weyl's tensor); in 7) the covectors
%$a$ and $b$ take specific form dependent on $R$ and derivatives.
%Conditions 1),3),4),7) imply that the Ricci tensor is $R$-compatible.\\
%The Codazzi condition has been extended and studied in the literature for 
%symmetric tensors $b_{ij}$. In this paper we extend 4) (to Quasi Codazzi
%tensor) and Sinyukov's equation (to generalized Weyl tensors). 

%The same approach can be used for a symmetric tensor $b_{jk}$. However,
%if $b$ is Ricci-like, i.e. $\nabla_mb^m{}_j=\frac{1}{2}\nabla_jb^m{}_m$, then 
%there is no difference with the preceding classification.
%
%\begin{prop}
For a simmetric tensor $b_{kl}$ with $\nabla_jb_{kl}\neq 0$, 
the tensor $\nabla_jb_{kl}$ can be decomposed into $O(n)$
invariant terms, where $\mathscr B^0_{jkl}$ is traceless 
($\mathscr B^0_{jk}{}^j=\mathscr B^0_{kj}{}^j=0$) \cite{Hamermesh, krupka}: 
\begin{gather}
\nabla_jb_{kl}= \mathscr B^0_{jkl}+ A_j g_{kl}+B_kg_{jl}+B_l g_{jk}\\
A_j=\frac{(n+1)\nabla_j b^m{}_m-2\nabla_mb^m{}_j}{n^2+n-2},\qquad
B_j=-\frac{\nabla_j b^m{}_m-n\nabla_mb^m{}_j}{n^2+n-2} 
\end{gather}
The traceless tensor can then be written as a sum of orthogonal components 
\cite{Lovelock}:
\begin{equation}
\mathscr B^0_{jkl}= \frac{1}{3}\left[\mathscr B^0_{jkl}+\mathscr B^0_{klj}+
\mathscr B^0_{ljk}\right]+ \frac{1}{3}\left [\mathscr B^0_{jkl}-
\mathscr B^0_{kjl}\right]+\frac{1}{3}\left [\mathscr B^0_{jlk}-
\mathscr B^0_{ljk}\right]
\end{equation}

The orthogonal subspaces classify the $O(n)$ invariant equations that are 
linear in $\nabla_jb_{kl}$. The trivial subspace: $\nabla_jb_{kl}=0$.
The subspace $\mathcal I$ (we follow Gray's notation, \cite{Gray}) 
where $\mathscr B^0_{jkl}=0$:
$$  \nabla_j b_{kl}= A_j g_{kl}+B_kg_{jl}+B_l g_{jk}.  $$
The complement $\mathcal I^\perp $ is characterized by $A_j,B_j=0$ i.e. 
$\nabla_j b_{kl}$ is traceless. This gives two invariant equations:
$\nabla_jb^j{}_l=0$, and $\nabla_jb^m{}_m=0$. Since $\nabla_jb_{kl}
= \mathscr B^0_{jkl}$, the structure of $\mathscr B^0$ 
specifies two orthogonal subspaces  
$\mathcal I^\perp =\mathcal A\oplus \mathcal B$. In $\mathcal A$:\\
$$   \nabla_jb_{kl}+\nabla_k b_{lj}+\nabla_lb_{jk}=0.  $$ 
In $\mathcal B$: 
$$ \nabla_jb_{kl}-\nabla_k b_{jl}=0. $$ 
%The conditions V and VI imply III and IV.\\
The subspace $\mathcal I\oplus \mathcal A$ contains tensors with traceless
part $\nabla_jb_{kl}-A_j g_{kl}-B_kg_{jl}-B_l g_{jk}$ that solves the cyclic 
condition:
 %$[\nabla_jb_{kl}+(A_j+2B_j)g_{kl}]+cyclic =0$ i.e.\\
$$ 
[\nabla_jb_{kl}-
\frac{1}{n+2}(\nabla_jb^m{}_m +2\nabla_mb^m{}_j)g_{kl}]+cyclic =0. 
$$
The subspace $\mathcal I\oplus \mathcal B$ contains tensors with traceless
part that solves the Codazzi condition:
%8) $[\nabla_jb_{kl}+(A_j-B_j)g_{kl}]- (j-k) =0$ i.e.\\ 
$$
[\nabla_jb_{kl}-\frac{1}{n-1}(\nabla_jb^m{}_m -\nabla_mb^m{}_j)g_{kl}]=
[\nabla_kb_{jl}-\frac{1}{n-1}(\nabla_kb^m{}_m -\nabla_mb^m{}_k)g_{jl}] $$

Accordingly, the Codazzi deviation tensor has the (unique) decomposition
in irreducible components 
\begin{equation}
\mathscr C_{jkl}= \mathscr C^0_{jkl}+\lambda_jg_{kl}  
-\lambda_k g_{jl} ,\qquad \lambda_j=A_j-B_j=
\frac{\nabla_j b^m{}_m-\nabla_mb^m{}_j}{n-1}
\label{tracelessC}
\end{equation}
where $\mathscr C^0$ is traceless. Eq.\eqref{CCCcomp} becomes
\begin{align}
 b_{im} R_{jkl}{}^m + b_{jm} R_{kil}{}^m + b_{km}R_{ijl}{}^m=
\nabla_i\mathscr C^0_{jkl}+\nabla_j\mathscr C^0_{kil}+\nabla_k\mathscr C^0_{ijl}
\label{nablaCprime}\\
+g_{il}(\nabla_j\lambda_k-\nabla_k\lambda_j)+
g_{jl}(\nabla_k\lambda_i-\nabla_i\lambda_k)+
g_{kl}(\nabla_i\lambda_j-\nabla_j\lambda_i)\nonumber
\end{align}
There are only two orthogonal invariant cases:\\
- $\mathscr C^0_{jkl}=0$, then $b$ is $R-$compatible if and only if $\lambda $
is closed. If $b$ is the Ricci tensor, this requirement gives Nearly
conformally symmetric $(NCS)_n$ manifolds, that were introduced by Roter
\cite{Roter}.\\
-$\nabla_j b^m{}_m-\nabla_mb^m{}_j=0$ then $b$ is $R-$compatible if and only
if $\mathscr C=\mathscr C^0$ solves the second Bianchi dentity. If $b$ is
the Ricci tensor, this corresponds to $\nabla_j R=0$.
%Therefore, if $\lambda $ is closed, $b_{ij}$ is $R$-compatible if and only if
%the traceless component $\mathscr C^0_{jkl}$ solves the second Bianchi 
%identity. Note that $\lambda $ is closed if and only if 
%$\nabla_i \nabla_m b_j{}^m= \nabla_j\nabla_m b_i{}^m$. 

\begin{remark}
The decomposition \eqref{tracelessC} for the deviation of the Ricci tensor
turns out to be
\begin{equation}
\mathscr C_{jkl}= -\frac{n-2}{n-3} \nabla_m C_{jkl}{}^m
+\frac{1}{2(n-1)}\left[g_{kl}\nabla_j R -g_{jl} \nabla_k R\right ]
\label{RicciWeyl}
\end{equation}
where $C_{jkl}{}^m$ is the conformal curvature tensor, or Weyl's tensor.
In this case the $\lambda $ covector is closed.
\end{remark}

\section{Riemann compatibility: general properties}
The existence of a Riemann compatible tensor has various implications. 
A first one is the existence of a generalized curvature 
tensor. This leads to the generalization of Derdzinski-Shen theorem and 
other relations that were obtained for Codazzi tensors.\\

We need the definition, from Kobayashi and Nomizu's book \cite{[10]}:
\begin{definition}\label{def2.1}
A tensor $K_{ijlm}$ is a {\em generalized curvature tensor}
(or, briefly, a $K$-tensor) if it has the symmetries of the Riemann curvature 
tensor:\\
a) $K_{ijkl} =-K_{jikl} =-K_{ijlk}$,\\ 
b) $K_{ijkl} =K_{klij}$,\\
c) $K_{ijkl} +K_{jkil} +K_{kijl} =0$ (first Bianchi identity).
\end{definition}
\noindent
It follows that the tensor $K_{jk}=: -K_{mjk}{}^m$ is symmetric.

\begin{thrm}
If $b$ is $R$-compatible then $K_{ijkl}=:R_{ijpq}b^p{}_kb^q{}_l$ is a 
$K$-tensor.
\begin{proof}
a) For example: 
$K_{ijlk}= R_{ijrs}b_l{}^r b_k{}^s = R_{ijsr} b_l{}^s b_k{}^r=
- R_{ijrs}b_l{}^s b_k{}^r =-K_{ijkl}$.
Property c) follows from \eqref{RC}: 
$K_{ijkl} +K_{jkil} +K_{kijl} = 
R_{ijrs}b_k{}^r b_l{}^s + R_{jkrs}  b_i{}^r b_l{}^s + R_{kirs} b_j{}^r b_l{}^s
=  ( R_{jis}{}^r b_{kr}+ R_{kjs}{}^r b_{ir}+ R_{iks}{}^r b_{jr}) b_l{}^s=0 $.
Property b) follows from c): $K_{ijkl} +K_{jkil} +K_{kijl} =0$. 
Sum the identity over cyclic permutations of all indices $i,j,k,l$ and use the 
symmetries a).\\ 
It is easy to see that a first Bianchi identity holds also for the last 
three indices: $K_{ijkl} +K_{iklj} +K_{iljk} =0$.
\end{proof}
\end{thrm}
\noindent

The next result remarks the relevance of the local basis of eigenvectors 
of the Ricci tensor. Another symmetric contraction of the Riemann tensor was 
introduced by Bourguignon \cite{Bou}:
\begin{equation} 
\text{\r{R}}_{ij}=: b^{pq}R_{pijq}.
\end{equation}
\begin{thrm} \label{Riccicomp}
If $b$ is $R$-compatible then:\\
1) $b_{im} R_j{}^m-b_{jm} R_i{}^m=0$,\\
2) $b_{im}\text{\r{R}}_j{}^m -  b_{jm}\text{\r{R}}_i{}^m =0$
\begin{proof}
The first identity is proven by transvecting \eqref{RC} with $g^{kl}$.
The second one is a restatement of the symmetry of the tensor $K_{ij}$.
\end{proof}
\end{thrm}
\noindent
\begin{remark} A) Identies 1 and 2 are here obtained directly from 
$R$-compatibility. Bourguignon \cite{Bou} obtained them from
Weitzenb\"ock's formula for Codazzi tensors, and Derdzinski and Shen \cite{[7]}
from their theorem.\\
B) As the symmetric matrices $b_{ij}$, $R_{ij}$, \r{R}$_{ij}$ commute, 
they share at each point of the manifold an orthonormal set of $n$ 
eigenvectors.\\
C) If $b'$ is a symmetric tensor that commutes with a Riemann compatible $b$, 
then it can be shown that \r{R}$'_{ij}=: b'^{pq}R_{pijq}$ commutes with $b$. 
\end{remark}

Finally, this Veblen-type identity holds:
\begin{prop} 
If $b$ is $R$-compatible, then:
\begin{equation} 
b_{im}R_{jlk}{}^m +b_{jm}R_{kil}{}^m +b_{km}R_{lji}{}^m +b_{lm}R_{ikj}{}^m = 0
\label{veblenb}
\end{equation}
\begin{proof}
Write four equations \eqref{RC} with cycled indices $i,j,k,l$ 
and sum up, and use the first Bianchi identity.
\end{proof}
\end{prop}

\section{Pure Riemann tensors and Pontryagin forms}
Riemann compatibility and nondegeneracy of the eigenvalues of $b$ imply
directly that the Riemann tensor is {\em pure} and Pontryagin forms vanish.\\
We quote two results from Maillot's paper \cite{Maillot}: 
\begin{definition}
In a Riemann manifold $\mathscr M_n$, the Riemann curvature 
tensor is pure if at each point of the manifold there is an orthonormal basis 
of $n$ tangent vectors $X(1),\ldots ,X(n)$, $X(a)^iX(b)_i = \delta_{ab}$, 
such that 
the tensors $X(a)^i\wedge X(b)^j=: X(a)^iX(b)^j-X(a)^jX(b)^i$, $a<b$, 
diagonalize it: 
\begin{equation}
R_{ij}{}^{lm} X(a)^i\wedge X(b)^j = \lambda_{ab}  X(a)^l\wedge X(b)^m 
\end{equation}  
\end{definition}
\begin{thrm}
If a Riemannian manifold has pure Riemann curvature tensor, 
then all Pontryagin forms vanish.
\end{thrm}

Consider the maps on tangent vectors, built with the Riemann tensor,
\begin{align*}
&\omega_4(X_1\ldots X_4)=R_{ija}{}^bR_{klb}{}^a(X_1^i\wedge X_2^j)
(X_3^k\wedge X_4^l), \\
&\omega_8(X_1\ldots X_8)=R_{ija}{}^bR_{klb}{}^c R_{mnc}{}^dR_{pqd}{}^a
(X_1^i\wedge X_2^j)\cdots(X_7^p\wedge X_8^q),\\
&\ldots\ldots
\end{align*}
They are antisymmetric under exchange of vectors in the single pairs, 
and for cyclic permutation of pairs. The {\em Pontryagin forms} \cite{Nakahara} 
$\Omega_{4k}$ result from total antisymmetrization of $\omega_{4k}$:
$ \Omega_{4k} (X_1\ldots X_{4k})=\sum_P (-1)^P \omega_{4k} 
(X_{i_1}\ldots X_{i_{4k}})$
where $P$ is the permutation taking $(1\ldots 4k)$ to $(i_1\ldots i_{4k})$.
$\Omega_{4k}=0$ if two vectors repeat, intermediate forms $\Omega_{4k-2}$ 
vanish identically.\\ 
Pontryagin forms on generic tangent vectors are linear combinations of 
forms evaluated on basis vectors.\\ 
If the Riemann tensor is pure, all Pontryagin forms on the basis 
of eigenvectors of the Riemann tensor vanish. For example,
if $X,Y,Z,W$ are orthogonal:
$\omega_4 (XYZW)= \lambda_{XY}\lambda_{ZW}(X^a\wedge Y^b)(Z_b\wedge W_a)=0 $ 
and $\Omega_4(XYZU)=0$. \\

A consequence of the extended Derdzinski-Shen theorem \ref{extDS} is the 
following: 

\begin{thrm}\label{distinct}
If a symmetric tensor field $b_{ij}$ exists, that is $R$-compatible 
and has distinct eigenvalues at each point of the manifold, then the Riemann
tensor is pure and all Pontryagin forms vanish. 
\begin{proof}
At each point of the manifold the symmetric matrix $b_{ij}(x)$ is 
diagonalized by $n$ tangent orthonormal vectors $X(a)$, with 
distinct eigenvalues. Since $b$ is $R$-compatible, theorem \ref{extDS} 
holds and, because of antisymmetry of $R$ in first two indices:
$$0=R_{ij}{}^{kl}X(a)^i\wedge X(b)^jX(c)_k, \quad  a\neq b\neq c. $$ 
This means that all column vectors of the matrix 
$V(a,b)^{kl}=R_{ij}{}^{kl}X(a)^i\wedge X(b)^j$ are orthogonal to vectors $X(c)$ 
i.e. they belong to the subspace spanned by $X(a)$ and $X(b)$. Because of 
antisymmetry in indices $k,l$, it is necessarily 
$V(a,b)=\lambda_{ab} X(a)\wedge X(b)$, i.e. the Riemann tensor is pure.  
\end{proof}
\end{thrm}

This property has been checked by Petersen \cite{Petersen} 
in various examples with 
rotationally invariant metrics, by giving explicit orthonormal frames 
such that $R(e_i,e_j)e_k=0$.

\subsection{Two and three dimensional manifolds}
Riemannian manifolds of dimension $n=2$ and $n=3$ are special, as the 
Riemann tensor is expressible in terms of the Ricci and metric tensors.
Therefore, Riemann-compatibility and ensuing pureness of the Riemann tensor 
can be established by simple means. \\
$n=2$) $R_{jklm}= R_{jl}g_{km}-g_{jm} R_{kj}$. Explicit evaluation 
proves that any symmetric tensor $b$ is Riemann compatible.\\
$n=3$) $R_{jklm}=g_{jl}R_{km}+g_{km}R_{jl}-g_{kl}R_{jm}-g_{jm}R_{kl}-\frac{R}{2}
(g_{jl}g_{km}-g_{jm}g_{kl})$. Then, for any symmetric tensor $b$ it is:
\begin{align*}
 b_{im} R_{jkl}{}^m+b_{jm} R_{kil}{}^m+b_{km} R_{ijl}{}^m= 
g_{kl}(b_{jm}R_i{}^m -b_{im}R_j{}^m)\\
+g_{il}(b_{km}R_j{}^m - b_{jm}R_k{}^m)+g_{jl}(b_{im}R_k{}^m - b_{km}R_i{}^m) 
\end{align*}
Thus in $n=3$ the Ricci tensor is always $R$-compatible. Moreover, if $b$ 
commutes with the Ricci tensor, then $b$ is $R$-compatible. 
Since a symmetric tensor that commutes with the Ricci tensor can always be 
constructed, with arbitrarily chosen distinct eigenvalues, by theorem 
\ref{distinct} we conclude:
\begin{prop}\label{2and3}
In Riemannian manifolds of dimension $n=2$ and $n=3$ the 
Riemann tensor is pure. 
\end{prop}

\subsection{Quasi-constant curvature spaces}
The same conclusions can be drawn  in any dimension $n$ for quasi-constant 
curvature spaces. They were introduced by Chen and 
Yano \cite{Chen_Yano} and are defined by a Riemann tensor with the
following structure: 
\begin{align}
R_{jklm} =p[g_{mj}g_{kl}-g_{mk}g_{jl}]+q[g_{mj}t_kt_l -g_{mk}t_jt_l+g_{kl}t_mt_j-
g_{jl}t_mt_k]
\end{align}
$p$ and $q$ are scalar functions. The first term describes constant 
curvature, the second one contains a vector field with $t_kt^k=1$.\\ 
The following identity holds:
\begin{align}
b_i{}^m R_{jklm} + b_j{}^m R_{kilm} + b_k{}^m R_{ijlm} =
q[g_{kl}(t_jb_i{}^mt_m-t_ib_j{}^mt_m)\label{Qcc} \\
+g_{il}(t_kb_j{}^mt_m-t_jb_k{}^mt_m)+g_{jl}(t_ib_k{}^mt_m-t_kb_i{}^mt_m)]
\nonumber
\end{align}
Contraction with $g^{kl}$ gives:
$ -b_i{}^m R_{jm} + b_j{}^m R_{im} = q(n-2)(t_jb_i{}^mt_m-t_ib_j{}^mt_m)$.
Therefore, if $b$ commutes with the Ricci tensor and $n\neq 2$, the r.h.s.
is zero and, by \eqref{Qcc}, $b$ is $R$-compatible. Then
the Riemann tensor is pure and all Pontryagin forms vanish.

\section{Structures for Riemann compatibility}
Some differential structures are presented that yield Riemann compatibility.
Of particular interest are geodesic mappings, which leave the condition for 
$R$-compatibility form-invariant, and generate $R$-compatible
tensors.
\subsection{Quasi Codazzi tensors} 
Let $b_{ij}$ be a symmetric tensor that solves the Codazzi condition deformed 
by a closed gauge field \cite{manticaDS}:
\begin{equation}\label{gauge}
(\nabla_j - \beta_j) b_{kl} = (\nabla_k - \beta_k) b_{jl} 
\end{equation}
The Codazzi deviation is  $\mathscr C_{jkl}= \beta_j b_{kl}- \beta_k b_{jl}$,
and $b$ is $R$-compatible.\\
Since $\beta_i =\nabla_i\xi$, the gauge field cancels for 
$b_{ij}=e^\xi b'_{ij}$, where $b'$ is a Codazzi tensor.

Of this type are Weakly $b$-symmetric manifolds, defined by the recurrency
\begin{equation}
\nabla_i b_{kl}= A_i b_{kl}+ B_k b_{il}+ D_l b_{ik} \label{wb}
\end{equation}
where $A$, $B$ and $D$ are covector fields. Eq.\eqref{gauge} is obtained 
for $\beta_i=A_i-B_i$, and $b$ is Riemann compatible if $A-B$ is closed.\\
Examples are: Weakly Ricci-symmetric manifolds, where $b_{ij}=R_{ij}$ 
\cite{[12], Weakly}, Weakly and pseudo Z-symmetric manifolds,
where $b_{ij}$ is a $Z$-tensor \cite{Weakly,Suh1}. Another example are 
manifolds with a recurrent generalized curvature tensor \cite{[12]}: 
$\nabla_iK_{jkl}{}^m= A_i K_{jkl}{}^m$, then $b_{kl}=: K_{kml}{}^m\neq 0$ 
has the form \eqref{wb}. 

\subsection{Pseudo-$K$symmetric manifolds}
They are characterized by a generalized curvature tensor $K$ such that 
(\cite{[4],closedness}) 
$$\nabla_iK_{jkl}{}^m = 2A_iK_{jkl}{}^m + A_j K_{ikl}{}^m + A_k K_{jil}{}^m+
A_l K_{jki}{}^m +A^m K_{jkli},$$
%which ensures validity of the second Bianchi identity. 
The tensor $b_{jk}=:K_{jmk}{}^m$ is symmetric. It is $R-$compatible if
its Codazzi deviation 
$\mathscr C_{ikl}= A_i b_{kl}-A_k b_{il} + 3A_m K_{ikl}{}^m $ 
fulfills the II Bianchi identity. This is ensured by $A_m$ being 
concircular, i.e. $\nabla_iA_m=A_iA_m+\gamma \,g_{im}$.  

\subsection{Generalized Weyl tensors}
A Riemannian manifold is a $(NCS)_n$ \cite{Roter} if the Ricci tensor 
satisfies $\nabla_jR_{kl}-\nabla_kR_{jl} = \frac{1}{2(n-1)}[g_{kl}\nabla_j\, R
-g_{jl}\nabla_k\,R  ]$. As such, the Ricci tensor is the Weyl tensor, 
and the left hand side is its Codazzi deviation. This condition, by
\eqref{RicciWeyl}, is equivalent to $\nabla_mC_{jkl}{}^m=0$.\\
This suggests a class of deviations of a symmetric
tensor $b$ with $\mathscr C^0_{jkl}=0$ in \eqref{tracelessC}:
\begin{equation}\label{genWeyl}
\mathscr C_{jkl}=\lambda_j g_{kl}-\lambda_k g_{jl} 
\end{equation}

\begin{prop}
$b$ is $R$-compatible if and only if $\lambda_i$ is closed.
\begin{proof}
Transvect \eqref{nablaCprime} with $g^{kl}$ and obtain:
$ -b_i{}^m R_{jm} + b_j{}^m R_{im} = 
(n-2)(\nabla_i\lambda_j-\nabla_j\lambda_i)$.
Then $b$ commutes with the Ricci tensor iff $\lambda $ is closed and, by 
the previous equation, $b$ is $R$-compatible. 
\end{proof}
\end{prop}

An example is provided by 
spaces with 
\begin{equation}\label{nablabg}
\nabla_j b_{kl}=A_j g_{kl}+B_k g_{jl}+B_l g_{jk},
\end{equation} 
where $\mathscr C_{jkl}= \lambda_jg_{kl}-\lambda_kg_{jl}$ with 
$\lambda_j=A_j-B_j$. Sinyukov manifolds \cite{Sinyukov} are of this sort, 
with $b_{ij}$ being the Ricci tensor itself.

\subsection{Geodesic mappings}
Riemann compatible tensors arise naturally in the study of geodesic mappings,
i.e. mappings that preserve geodesic lines. Their importance arise from the 
fact that Sinyukov manifolds are $(NCS)_n$ manifolds and they always admit a 
nontrivial geodesic mapping.\\
Geodesic mappings preserve Weyl's 
projective curvature tensor \cite{Sinyukov}. 
We show that they also preserve the form of the compatibility relation.

A map $f:\; (\mathscr M_n,g)\to (\mathscr M_n,\overline g)$  is 
{\em geodesic} if and only if Christoffel symbols are linked by 
$\overline\Gamma_{ij}^k = \Gamma_{ij}^k +\delta_i^k X_j + \delta_j^k X_i$
%\label{geoGamma}\end{equation}
where, on a Riemannian manifold, $X$ is closed ($\nabla_iX_j= \nabla_jX_i$).
The condition %\eqref{geoGamma} 
is equivalent to:
\begin{equation}
\nabla_k \overline g_{jl} = 2X_k\overline g_{jl} + X_j\overline g_{kl}+
X_l\overline g_{kj}\label{geog} 
\end{equation}
which has the form \eqref{nablabg}. The corresponding relation among 
Riemann tensors is
\begin{equation}\label{Riemgeod}
\overline R_{jkl}{}^m = R_{jkl}{}^m  +\delta_j^m P_{kl} - \delta_k^m P_{jl}
\end{equation} 
where $P_{kl}=\nabla_kX_l-X_kX_l$ is the {\em deformation} tensor. The symmetry
$P_{kl}=P_{lk}$ is ensured by closedness of $X$.
\begin{prop}
Geodesic mappings preserve $R$-compatibility
\begin{equation}
b_{im} \overline R_{jkl}{}^m+b_{jm} \overline R_{kil}{}^m+
b_{km} \overline R_{ijl}{}^m =  
b_{im} R_{jkl}{}^m+b_{jm} R_{kil}{}^m+b_{km} R_{ijl}{}^m 
\end{equation}
where $b$ is a symmetric tensor.
\begin{proof}
Let's show that the difference of the two sides is zero. Eq.(\ref{Riemgeod})
gives:\\ 
$b_{im} (\delta_j^m P_{kl} - \delta_k^m P_{jl}) +b_{jm} 
(\delta_k^m P_{il} - \delta_i^m P_{kl})+b_{km} 
(\delta_i^m P_{jl} - \delta_j^m P_{il})$\\
$=b_{ij}P_{kl} - b_{ik} P_{jl} +b_{jk} P_{il} - b_{ji}P_{kl}+
b_{ki} P_{jl} - b_{kj} P_{il}=0$
\end{proof}
\end{prop}
\noindent
Since $\overline g$ is trivially $\overline R$-compatible (first Bianchi
identity), form invariance implies:
\begin{cor}
$\overline g $ is $R$-compatible.
\end{cor}

Sinyukov \cite{Sinyukov} (see also \cite{Mikes,Formella}) showed 
that a manifold admits a geodesic mapping if and only if 
there are a scalar field $\varphi$ and a symmetric non singular tensor $b_{ij}$ 
such that:  
$$ \nabla_k b_{jl} = g_{kl} \nabla_j\varphi  + g_{kj}\nabla_l \varphi. $$
%The conditions for $b$ to have the given form are:
%\begin{align}
%n \nabla_r\lambda_j =\mu g_{jr} -R_r{}^m b_{mj}-R^i{}_{rj}{}^m b_{im}\\
%-2(n+1)\lambda^m R_{km}=(n-1)\nabla_k\mu + 
%b^{im}(2\nabla_iR_{km}-\nabla_k R_{im})
%\end{align}
%with $\lambda_j=\nabla_j\varphi$.
%
The Codazzi deviation of $b$, 
$\mathscr C_{jkl}= g_{kl}\nabla_j\varphi -g_{jl}\nabla_k\varphi $,
has the form \eqref{genWeyl}. Therefore $b$ is $R$-compatible.\\

\section{Generalized curvature tensors.}
Several results that are valid for the Riemann tensor with a
Riemann compatible tensor, extend to generalized curvature tensors
$K_{ijkl}$ (hereafter referred to as $K$-tensors) with a $K$-compatible 
symmetric tensor $b_{jk}$. The classical curvature tensors are $K$-tensors.
The compatibility with the Ricci tensor is then examined.   
\begin{definition}\label{def2.2}
A symmetric tensor $b_{ij}$ is {\em K-compatible} if
\begin{equation}
b_{im} K_{jkl}{}^m+b_{jm} K_{kil}{}^m+b_{km} K_{ijl}{}^m=0.\label{Kcom}
\end{equation}
\end{definition}
\noindent
The metric tensor is always $K$-compatible, 
as \eqref{Kcom} then coincides with the first Bianchi identity for $K$.

\begin{prop}\label{lem2.2}
If $K_{ijlm}$ is a $K$-tensor and $b_{kl}$ is $K$-compatible, 
then $\hat K_{ijkl} =: K_{ijrs} b_k{}^r b_l{}^s$ is a $K$-tensor.
\end{prop}
We quote without proof the extension of Derdzinski and Shen theorem for
generalized curvature tensors \cite{manticaDS}:
\begin{thrm}
Suppose that $K_{ijkl}$ is a $K$-tensor, and a 
symmetric $K$-compatible tensor $b_{ij}$ exists.
Then, if $X$, $Y$ and $Z$ are three eigenvectors of the matrix $b_r{}^s$
at a point $x$ of the manifold, with eigenvalues $\lambda $, $\mu $ and $\nu $,
it is $X^i \, Y^j\, Z^k \, K_{ijkl} =0$ 
%\begin{equation}
%X^i \, Y^j\, Z^k \, K_{ijkl} =0
%\end{equation}
provided that both $\lambda $ and $\mu $ are different from $\nu$. 
\end{thrm}

\begin{prop}\label{Kh}
If $b$ is $K-$ compatible, and $b$ commutes with a tensor $h$, 
then the symmetric tensor {\r K}$_{kl}=: K_{jklm}h^{jm}$ commutes with 
$b$.
\begin{proof}
Multiply relation of $K$ compatibility for $b$ by $h^{kl}$. The last term
vanishes for symmetry. The remaining terms give the null commutation relation. 
\end{proof}
\end{prop}

In ref.\cite{[12]} (prop.2.4) we proved that a generalization of 
Lovelock's identity \eqref{lovelock} holds for certain $K$-tensors,
that include all classical curvature tensors:
\begin{prop}\label{thrm1.3}
Let $K_{jkl}{}^m$ be a $K$-tensor with the property
\begin{equation}
\nabla_m K_{jkl}{}^m = \alpha \,\nabla_m R_{jkl}{}^m+
\beta \, \left (a_{kl}\nabla_j - a_{jl}\nabla_k\right ) \varphi ,\label{eq1.9}
\end{equation}
where $\alpha $, $\beta $ are non zero constants, $\varphi $ is a real 
scalar function and $a_{kl}$ is a Codazzi tensor. Then:
\begin{eqnarray}
&&\nabla_i\nabla_m K_{jkl}{}^m+ \nabla_j\nabla_m K_{kil}{}^m +\nabla_k\nabla_m 
K_{ijl}{}^m\label{eq1.10} \\
&&\qquad =-\alpha (R_{im} R_{jkl}{}^m+R_{jm} R_{kil}{}^m + R_{km} R_{ijl}{}^m).
\nonumber
\end{eqnarray}
\end{prop}
%(it also follows from eq. 4.8 in \cite{Bou})

\subsection{ABC curvature tensors}
A class of $K$-tensors with the structure \eqref{eq1.9} are the $ABC$ 
curvature tensors. They are combinations of the Riemann tensor and its 
contractions ($A$, $B$, $C$ are constants unless otherwise stated):
\begin{eqnarray}
K_{jkl}{}^m = R_{jkl}{}^m \,+\, A(\delta_j{}^m R_{kl}-\delta_k{}^m R_{jl}) 
+ B(R_j{}^m g_{kl}-R_k{}^m g_{jl})\\
+ C(R\delta_j{}^m g_{kl}-R\delta_k{}^m g_{jl})\nonumber
\end{eqnarray}
The canonical curvature tensors are of this sort:
\begin{itemize}
\item{{\em Conformal} tensor $C_{ijkl}$: $A=B=\frac{1}{n-2}$, 
$C= -\frac{1}{(n-1)(n-2)}$;}
\item{{\em Conharmonic} tensor $N_{ijkl}$: $A=B=\frac{1}{n-2}$, $C=0$;}
\item{{\em Projective} tensor: $P_{ijkl}$: $A= \frac{1}{n-1}$, $B=C=0$;}
\item{{\em Concircular} tensor: 
$\tilde C_{ijkl}$: $A-B=0$, $C= \frac{1}{n(n-1)}$.}
\end{itemize}
\begin{prop}
Let $K_{ijkl}$ be an $ABC$ tensor ($A$, $B$, $C$ may be scalar
functions) and $b_{ij}$ a symmetric tensor; \\
1) if $b$ is $R$-compatible then $b$ is $K$-compatible.\\
2) if $b$ is $K$-compatible and $B\neq \frac{1}{n-2}$ then $b$ is
$R$-compatible.
\begin{proof}
The following identity holds for $ABC$ tensors and a symmetric tensor $b$:
\begin{align}
b_{im} K_{jkl}{}^m + b_{jm} K_{kil}{}^m + b_{km} K_{ijl}{}^m 
= b_{im} R_{jkl}{}^m + b_{jm} R_{kil}{}^m + b_{km} R_{ijl}{}^m \label{RK}\\
+ B\left [ 
g_{kl} (b_{im} R_j{}^m - b_{jm} R_i{}^m ) + g_{il} (b_{jm} R_k{}^m - b_{km} 
R_j{}^m ) + g_{jl} (b_{km} R_i{}^m - b_{im} R_k{}^m )\right ].\nonumber
\end{align}
1) by theorem \ref{Riccicomp}, if $b$ is $R$-compatible then it commutes with 
the Ricci tensor, and $K$-compatibility follows.\\
2) if $b$ is $K$-compatible it commutes with $K_{ij}$. Contraction with 
$g^{kl}$ gives:
\begin{equation*}
b_{im}K_j{}^m - b_{jm} K_i{}^m = (b_{im}R_j{}^m - b_{jm} R_i{}^m)[1-B(n-2)],
\end{equation*}
then, if $B\neq \frac{1}{n-2}$, $b$ commutes with the Ricci tensor and  
by \eqref{RK} it is $R$-compatible.
\end{proof}
\end{prop}

%\begin{cor}
%If the Ricci tensor is $R$-compatible, then it is $ABC$-compatible.
%\end{cor}

%\begin{cor}
%If $K$ is a flat $ABC$ tensor, then a symmetric tensor $b_{ij}$ is
%$R$-compatible iff it commutes with the Ricci tensor. If this occurs,
%the $ABC$ tensor is pure and Pontryagin forms vanish.
%\end{cor}

\begin{prop}\label{prop66}
Let $K$ be an $ABC$ tensor with constant $A\neq 1$ and $B$. If
 \begin{eqnarray}
\nabla_i\nabla_m K_{jkl}{}^m+ \nabla_j\nabla_m K_{kil}{}^m +\nabla_k\nabla_m 
K_{ijl}{}^m=0 
\end{eqnarray}
then the Ricci tensor is $K$-compatible.
\begin{proof}
If $A$ and $B$ are constants, one evaluates
\begin{eqnarray}\label{divABC}
\nabla_m K_{jkl}{}^m = (1-A) \nabla_m R_{jkl}{}^m + \frac{1}{2}(B +2C)\left(
 g_{kl}\nabla_jR -  g_{jl}\nabla_k R \right ),
\end{eqnarray}
Lovelock's identity \eqref{lovelock} for the Riemann tensor implies
\begin{eqnarray}
&&\nabla_i\nabla_m K_{jkl}{}^m+ \nabla_j\nabla_m K_{kil}{}^m +\nabla_k\nabla_m 
K_{ijl}{}^m \label{abcxx}\\
&&\qquad =-(1-A)(R_{im} R_{jkl}{}^m+R_{jm} R_{kil}{}^m + R_{km} R_{ijl}{}^m).
\nonumber
\end{eqnarray}
In the r.h.s. the Riemann tensor can be replaced by tensor $K$ by \eqref{RK}
written for the Ricci tensor.
\end{proof}
\end{prop}

Sufficient conditions are: $K$ is harmonic, $K$ is recurrent 
(with closed recurrency parameter, see eq.(3.13) in \cite{[12]}). Note that 
prop. \ref{prop66} 
remains valid for the Weyl's conformal tensor, which is traceless.\\

%Given a non trivial Codazzi tensor $a_{kl}$, the $K$-tensor $K_{jkl}{}^m= 
%A\, R_{jkl}{}^m + B\,\varphi \,(\delta_j{}^m a_{kl} -\delta_k{}^m 
%a_{jl})$ 
%%.\label{eq1.11}\end{equation}
%is not $ABC$-type but satisfies \eqref{eq1.9} 
%\cite{Weakly,closedness}.
%The trace defines a $Z$ tensor 
%%\begin{equation}
%$Z_{kl} = -K_{mkl}{}^m = AR_{kl} + B(1-n)\varphi\, b_{kl} $
%%.\label{eq1.12}\end{equation}
%The notion of such $Z$ tensor arises naturally from the invariance property 
%of Lovelock's identity.

\section{Weyl-compatibility and General Relativity} 
In general relativity, the Ricci tensor is related to
the energy-momentum tensor by the Einstein equation: 
$R_{jl}=\frac{1}{2}Rg_{jl}+kT_{jl}$ with curvature $R=-2kT/(n-2)$ 
($T=T^k{}_k$).\\
The contracted II Bianchi identity gives
$$\nabla_mR_{jkl}{}^m= k \left(\nabla_k T_{jl}-\nabla_j T_{kl}\right )
+\frac{1}{2} \left( g_{jl} \nabla_k R -   g_{kl}\nabla_j R \right ). $$
Let $K$ be an $ABC$ tensor, with constant $A$, $B$, $C$. Its divergence 
\eqref{divABC} can be expressed in terms of the gradients of the
energy momentum tensor $T_{ij}$.
%\begin{eqnarray}
%\nabla_m K_{jkl}{}^m = k (1-A) \left(\nabla_k T_{jl}-\nabla_j T_{kl}\right )
%+\frac{k}{2-n}(1-A-B-2C)\left[g_{jl} \nabla_k T  - g_{kl}\nabla_j T \right ]
%\end{eqnarray}
In the same way Einstein's 
equations and \eqref{abcxx} give an equation which is local in the energy 
momentum tensor:
\begin{eqnarray}
&&\nabla_i\nabla_m K_{jkl}{}^m+ \nabla_j\nabla_m K_{kil}{}^m +\nabla_k\nabla_m 
K_{ijl}{}^m\label{eq1.10} \\
&&\qquad =-(1-A)k\left (T_{im} K_{jkl}{}^m+T_{jm} K_{kil}{}^m + T_{km} K_{ijl}{}^m
\right ).\nonumber
\end{eqnarray}
The Weyl tensor $C_{jkl}{}^m$ is the traceless part of the Riemann tensor, 
and it is an $ABC$ tensor. There are advantages in discussing General 
Relativity by taking the Weyl tensor as the fundamental geometrical 
quantity \cite{Berts, HawkingEllis, DeFelice}. The first equation
\eqref{divABC}
$$ \nabla_m C_{jkl}{}^m = 
k \frac{n-3}{n-2} \left[\nabla_k T_{jl}-\nabla_j T_{kl} +
\frac{1}{n-1}\left( g_{jl} \nabla_k T  - g_{kl}\nabla_j T \right )\right] $$
is reported in textbooks, as De Felice \cite{DeFelice},  
Hawking Ellis \cite{HawkingEllis}, Stephani \cite{Stephani}, and in the paper
\cite{Berts}.
Instead, a further derivation yields a Bianchi-like equation for the 
divergence, Eq.\eqref{eq1.10}, which contains 
no derivatives of the sources
\begin{eqnarray}
&&\nabla_i\nabla_m C_{jkl}{}^m+ \nabla_j\nabla_m C_{kil}{}^m +\nabla_k\nabla_m 
C_{ijl}{}^m \\
&&\qquad =-k\frac{n-3}{n-2} \left ( T_{im} C_{jkl}{}^m+T_{jm} C_{kil}{}^m + 
T_{km} C_{ijl}{}^m \right ).\nonumber
\end{eqnarray}
It can be viewed as a condition for Weyl-compatibility for the energy
momentum tensor.\\
In view of prop.\ref{Kh} and the previous equation, the following holds:
\begin{prop}
If $T_{ij}$ is Weyl-compatible, the symmetric tensor 
$\text{\r {C}}_{kl}=: T^{jm}C_{jklm}$ commutes with $T_{ij}$.
\end{prop}

In 4 dimensions, given a time-like velocity field $u^i$, Weyl's
tensor is projected in longitudinal (electric)
and transverse (magnetic) tensorial components \cite{Berts}
$$ E_{kl}=u^ju^m C_{jklm}, \qquad H_{kl}=\frac{1}{4}u^ju^m (\epsilon_{pqjk}
C^{pq}{}_{lm} + \epsilon_{pqjl} C^{pq}{}_{km}) $$
that solve equations that resemble Maxwell's equations with source.
Therefore, the tensor $E_{kl}=\text{\r {C}}_{kl}$ can be viewed as a 
generalized electric field. It coincides with
the standard definition if $T_{ij}=(p+\rho)u_iu_j+pg_{ij}$ (perfect fluid).
The generalized magnetic field is $H_{kl} = 
\frac{1}{4}T^{jm} (\epsilon_{pqjk}
C^{pq}{}_{lm} + \epsilon_{pqjl} C^{pq}{}_{km}) $.

%%%%%The subject is currently under investigation.

\begin{prop}
If $T_{kl}$ is Weyl compatible then $H_{kl}=0$. 
\begin{proof} From the condition for Weyl compatibility 
%$T^{im}C^{jk}{}_{lm}+T^{jm}C^{ki}{}_{lm}+T^{km}C^{ij}{}_{lm}=0$
we obtain 
$\epsilon_{ijkp}[T^{im}C^{jk}{}_{lm}+T^{jm}C^{ki}{}_{lm}+T^{km}C^{ij}{}_{lm}]=0 $.
The first and the second term are modified as follows:
\begin{align*}
&\epsilon_{ijkp} T^{im}C^{jk}{}_{lm}=\epsilon_{kijp}T^{km}C^{ij}{}_{lm}=
\epsilon_{ijkp}T^{km}C^{ij}{}_{lm}\\
&\epsilon_{ijkp}T^{jm}C^{ki}{}_{lm}=\epsilon_{jkip}T^{km}C^{ij}{}_{lm}=
\epsilon_{ijkp}T^{km}C^{ij}{}_{lm}.
\end{align*}
Then, since the sum becomes $\epsilon_{ijkp}T^{km}C^{ij}{}_{lm}=0$, then
the magnetic part of Weyl's tensor is zero.   
\end{proof}
\end{prop}

\end{document}